# GLOBAL EXPONENTIAL SAMPLED-DATA OBSERVERS FOR NONLINEAR SYSTEMS WITH DELAYED MEASUREMENTS

**Tarek Ahmed-Ali*, Iasson Karafyllis** and Francoise Lamnabhi-Lagarrigue*****

***Laboratoire GREYC CNRS-ENSICAEN**
**06 Boulevard du Marechal Juin 14050 Caen Cedex**
**(email: tarek.ahmed-ali@ensicaen.fr )**

****Department of Environmental Engineering, Technical University of Crete,**
**73100, Chania, Greece (email: ikarafyl@enveng.tuc.gr )**

*****Centre National de la Recherche Scientifique, CNRS-EECI SUPELEC,**
**3 rue Joliot Curie, 91192, Gif-sur-Yvette, France**
**(email: lamnabhi@lss.supelec.fr)**

## Abstract

This paper presents new results concerning the observer design for wide classes of nonlinear systems with both sampled and delayed measurements. By using a small gain approach we provide sufficient conditions, which involve both the delay and the sampling period, ensuring exponential convergence of the observer system error. The proposed observer is robust with respect to measurement errors and perturbations of the sampling schedule. Moreover, new results on the robust global exponential state predictor design problem are provided, for wide classes of nonlinear systems.

## 1. Introduction

In the last decades, the design of nonlinear observers for continuous systems with communication constraints has received a great attention. This interest is motivated by many engineering applications, such as sampled-data systems, network control systems (NCSs), and quantized systems. In the case of sampled-data systems, the output is available only at sampling instants. For linear systems it is usually possible to compute the discrete time model of the continuous time system. This is not the case for nonlinear systems where the exact discrete time model is generally not available.

In the nonlinear case there are several approaches in the literature for the design of sampled-data observers:

1) One approach is the design of a discrete observer by using a consistent approximation of the exact discretized model. This approach usually results to semi-global practical stability of the observation error. More details on this method can be found in [5] (see also references therein).
2) Another approach is based on a mixed continuous and discrete design. This approach has been inspired by Jazwinski in [13], who introduced the continuous discrete Kalman filter to solve a filtering problem for stochastic continuous-discrete time systems. In [8] the authors use this approach to write a discrete-continuous version of the well known high gain observer (see [9]). In [21], observers for a MIMO class of state affine systems where the dynamical matrix depends on the inputs have been designed when the inputs are regularly persistent. This work was extended to adaptive observers in [2]. In [6], a similar method has been used for a larger class of systems and applied to the observation of an emulsion copolymerization process. The observation of a class of systems with output injection has been treated in [22] and in [11] a high gain continuous-discrete observer has been developed by using constant observation gains. In [3], the authors extend the work of [11] to the discrete-time measurements case.
3) Recently, a new hybrid observer which uses an inter-sample predictor has been proposed in [14]. Sufficient conditions involving the sampling period have been derived by using a small gain approach. The results of [14] have been extended in [16].

On the other hand, the observation of systems with output delayed measurements has been considered in [10]. The authors proposed cascade predictors for a wide class of nonlinear systems to handle the delay of the output. Another kind of cascade predictors have been proposed in [7], [1]. The convergence has been derived by using *Lyapunov Krasovskii* tools. The design of observers for linear detectable systems with sampled and delayed measurements was treated in [12] by using a descriptor system approach and a Lyapunov Krasovskii functional. The authors proposed a hybrid observer, without using an inter-sample predictor, for a class of linear systems and derive sufficient conditions based on linear matrix inequalities to guarantee exponential convergence of the observation error. This idea has also been used in ([23], [25]) for some classes of nonlinear systems with sampled measurements.

In the work, we will present several results concerning the design of predictors and observers for certain classes of nonlinear systems with sampled and delayed measurements by using small gain arguments. We focus on nonlinear forward complete systems of the form:

$$\dot{x} = f(x,u), x \in \Re^n, u \in U \tag{1.1}$$

where $U \subseteq \Re^m$ is a non-empty set, $f : \Re^n \to \Re^n$ is a smooth vector field and the output is given by

$$y = h(x) \tag{1.2}$$

where $h : \Re^n \to \Re^k$ is a smooth mapping.

More specifically, Section 2 of the present paper provides a general result (Theorem 2.3), which guarantees that the combination of a sampled-data observer and a state predictor will yield a robust global exponential observer with sampled and delayed measurements. Robustness with respect to measurement errors and perturbations of the sampling schedule is guaranteed by Theorem 2.3. Section 3 of the present work focuses on globally Lipschitz systems, for which robust global exponential state predictors can be designed for arbitrary prediction horizon (Theorem 3.1). The linear time invariant case is treated as a special case of globally Lipschitz systems. It has to be noticed that the classes of observers that can be used in the proposed observer design include several well-known observers such as the high gain observers in [9] and the nonlinear observers designed in [24]. Section 4 of the present paper is devoted to nonlinear systems with a robustly globally asymptotically stable set, for which robust global exponential state predictors can be designed for sufficiently small prediction horizon (Proposition 4.1). The Appendix contains the proofs of the existence of exponential state predictors for the aforementioned classes of nonlinear systems.

**Notation.** Throughout this paper, we adopt the following notation:

* $\Re_+ := [0,+\infty)$. A partition of $\Re_+$ is a set $\pi = \{\tau_i\}_{i=0}^{\infty}$ with $\tau_0 = 0$, $\tau_{i+1} > \tau_i$ for all $i \ge 0$ and $\lim_{i \to +\infty} \tau_i = +\infty$.
* By $L^\infty(I;U)$, where $I \subseteq \Re$ is an interval and $U \subseteq \Re^m$ is a non-empty set, we denote the space of Lebesgue measurable and essentially bounded functions $u : I \to U$. For $u \in L^\infty(I;U)$ we denote by $\|u\|$ the essential supremum of $u$ on $I \subseteq \Re$. If $I \subseteq \Re$ is an unbounded interval then $L^\infty_{loc}(I;U)$ denotes the space of Lebesgue measurable and locally essentially bounded functions $u : I \to U$. By $C^0(A;U)$, where $A \subseteq \Re^n$ and $U \subseteq \Re^m$ are non-empty sets, we denote the space of continuous functions $u : A \to U$.
* By $[x]$ we denote the integer part of $x \in \Re$.
* For a vector $x \in \Re^n$, we denote by $x'$ its transpose and by $|x|$ its Euclidean norm. $A' \in \Re^{n \times m}$ denotes the transpose of the matrix $A \in \Re^{m \times n}$ and $|A|$ denotes the induced norm of the matrix $A \in \Re^{m \times n}$, i.e., $|A| = \sup\{|Ax| : x \in \Re^m, |x| = 1\}$.
* For a smooth vector field $\Psi : \Re^l \to \Re^n$, $D\Psi(z)$ denotes the Jacobian matrix of $\Psi$ at $z \in \Re^l$.
* For a function $V \in C^1(A;\Re)$, the gradient of $V$ at $x \in A \subseteq \Re^n$, denoted by $\nabla V(x)$, is the row vector $\nabla V(x) = \left[ \frac{\partial V}{\partial x_1}(x) \;\ldots\; \frac{\partial V}{\partial x_n}(x) \right]$. The Lie derivative of $V$ at $x \in A \subseteq \Re^n$ along the smooth vector field $f(x,u)$ is denoted by $L_f V(x,u) = \nabla V(x) f(x,u)$.

## 2. Robust global exponential observer with sampled and delayed measurements

Consider the forward complete system (1.1), (1.2). We first provide the definitions of the Robust Global Exponential Observer for (1.1) and the robust global exponential $r-$ predictor for (1.1).

**Definition 2.1:** *Consider the system*

$$\begin{aligned}&\dot{z} = F(z, y+v, u)\\&\hat{x} = \Psi(z)\\&z \in \Re^l, y \in \Re^k, v \in \Re^k, u \in U \subseteq \Re^m, \hat{x} \in \Re^n\end{aligned} \tag{2.1}$$

*where* $F : \Re^l \times \Re^k \times U \to \Re^l$ *and* $\Psi : \Re^l \to \Re^n$ *are smooth vector fields. We say that system (2.1) is a Robust Global Exponential Observer for (1.1), if there exist a non-decreasing function* $M : \Re_+ \to \Re_+$ *and constants* $\sigma > 0$, $\gamma \geq 0$ *such that for every* $(x_0, z_0, u, v) \in \Re^n \times \Re^l \times L^\infty(\Re_+; U) \times L^\infty_{loc}(\Re_+; \Re^k)$ *the solution* $(x(t), z(t))$ *of (1.1), (1.2) and (2.1) with initial condition* $(x(0), z(0)) = (x_0, z_0)$ *corresponding to inputs* $(u, v) \in L^\infty(\Re_+; U) \times L^\infty_{loc}(\Re_+; \Re^k)$ *exists for all* $t \geq 0$ *and satisfies the following estimate for all* $t \geq 0$:

$$|\hat{x}(t) - x(t)| \leq e^{-\sigma t} M(|x_0| + |z_0| + \|u\|) + \gamma \sup_{0\leq s\leq t}\left(e^{-\sigma(t-s)}|v(s)|\right) \tag{2.2}$$

**Definition 2.2:** *Let* $A(\Re_+; \Re^n)$ *be a non-empty subset of the functions* $z : \Re_+ \to \Re^n$ *which are absolutely continuous on every bounded interval of* $\Re_+$. *Consider a deterministic system of the form*

$$\begin{aligned}&\dot{\xi}(t) = F_p(\xi_t, u_t, z(t), \dot{z}(t))\\&\tilde{x}(t) = G(\xi_t, u_t, z(t))\\&\xi(t) \in \Re^q, \tilde{x}(t) \in \Re^n, u(t) \in U \subseteq \Re^m, z(t) \in \Re^n\end{aligned} \tag{2.3}$$

*where* $(\xi_t)(\theta) = \xi(t+\theta)$, $(u_t)(\theta) = u(t+\theta)$, *for* $\theta \in [-r, 0]$, $G : C^0([-r,0];\Re^q) \times L^\infty([-r,0];U) \times \Re^n \to \Re^n$ *is a continuous mapping and* $r > 0$ *is a positive constant. Assume that the mapping* $F_p : C^0([-r,0];\Re^q) \times L^\infty([-r,0];U) \times \Re^n \times \Re^n \to \Re^q$ *is continuous and such that for every* $(x_0, \xi_0, u, z) \in C^0([-r,0];\Re^n) \times C^0([-r,0];\Re^q) \times L^\infty([-r,+\infty);U) \times A(\Re_+;\Re^n)$ *the solution* $(x(t), \xi(t)) \in \Re^n \times \Re^k$ *of (1.1) and (2.3) with initial condition* $\xi(\theta) = (\xi_0)(\theta)$, $x(\theta) = (x_0)(\theta)$, $\theta \in [-r, 0]$ *and corresponding to inputs* $(u, z) \in L^\infty([-r,+\infty);U) \times A(\Re_+;\Re^n)$ *is unique, defined for all* $t \geq 0$ *and satisfies*

$$|\tilde{x}(t) - x(t)| \leq e^{-\sigma t} a\left(\|x_0\| + \|\xi_0\| + \|u\| + |z(0)| + \sup_{0\leq s\leq r}|\dot{z}(s)|\right) + P \sup_{0\leq s\leq t}\left(e^{-\sigma(t-s)}|z(s) - x(s-r)|\right),\ \forall t \geq 0 \tag{2.4}$$

*for certain non-decreasing function* $a : \Re_+ \to \Re_+$ *and certain constants* $P \geq 0$, $\sigma > 0$. *Then system (2.3) is called a robust global exponential* $r-$ *predictor for (1.1) with input* $z \in A(\Re_+; \Re^n)$.

We next provide the key hypotheses of this section. The hypotheses should be compared with the hypotheses of the main result in [14]. The hypotheses introduced here are more demanding but this is expected because here we consider systems with inputs, with delayed measurements and we require exponential convergence.

**Hypothesis (H1):** *System (1.1) admits a Robust Global Exponential Observer given by (2.1). Moreover, the system*

$$\begin{aligned}&\dot{z}(t) = F(z(t), w(t), u(t-r))\\&\dot{w}(t) = L_f h(\Psi(z(t)), u(t-r))\end{aligned} \tag{2.5}$$

*is forward complete for inputs* $u \in L^\infty([-r,+\infty);U)$. *Furthermore, the system*

$$\dot{z}(t) = F(z(t), w(t), u(t-r)) \tag{2.6}$$

*is forward complete for inputs* $(u,w) \in L^{\infty}([-r,+\infty);U) \times L^{\infty}_{loc}([-r,+\infty);\Re^{k})$.

**Hypothesis (H2):** *There exists a non-empty subset of the functions* $z:\Re_{+} \to \Re^{n}$ *which are absolutely continuous on every bounded interval of* $\Re_{+}$ *denoted by* $A(\Re_{+};\Re^{n})$ *such that system (1.1) admits a robust global exponential* $r-$*predictor for (1.1) with input* $z \in A(\Re_{+};\Re^{n})$, *which is given by (2.3). Moreover, for every* $(z_0,u,w) \in \Re^{l} \times L^{\infty}([-r,+\infty);U) \times L^{\infty}_{loc}([-r,+\infty);\Re^{k})$, *the output signal* $\hat{x}(t) = \Psi(z(t))$ *produced by the unique solution of (2.6) with initial condition* $z(0) = z_0$ *and corresponding to inputs* $(u,w) \in L^{\infty}([-r,+\infty);U) \times L^{\infty}_{loc}([-r,+\infty);\Re^{k})$ *is a function of class* $A(\Re_{+};\Re^{n})$.

**Hypothesis (H3):** *There exist a constant* $C > 0$, *a continuous function* $T:\Re^{n} \times \Re^{l} \to \Re_{+}$ *and a non-decreasing function* $N:\Re_{+} \to \Re_{+}$ *such that for every* $(x_0,z_0,u,v) \in \Re^{n} \times \Re^{l} \times L^{\infty}(\Re_{+};U) \times L^{\infty}_{loc}(\Re_{+};\Re^{k})$ *the solution* $(x(t),z(t))$ *of (1.1), (1.2) and (2.1) with initial condition* $(x(0),z(0)) = (x_0,z_0)$ *corresponding to inputs* $(u,v) \in L^{\infty}(\Re_{+};U) \times L^{\infty}_{loc}(\Re_{+};\Re^{k})$ *satisfies the following estimate for all* $t \geq T(x_0,z_0)$:

$$\left| L_f h(\hat{x}(t),u(t)) - L_f h(x(t),u(t)) \right| \leq e^{-\sigma t} N\left( |x_0| + |z_0| + \|u\| \right) + C \sup_{0 \leq s \leq t} \left( e^{-\sigma(t-s)} |v(s)| \right) \tag{2.7}$$

The main result of the section follows. The following theorem guarantees that there exists a global exponential sampled-data observer for system (1.1) under hypotheses (H1-3). Moreover, the observer is robust to measurement errors and perturbations of the sampling schedule.

**Theorem 2.3:** *Consider system (1.1) under hypotheses (H1-3). Let* $0 < b \leq B$ *be (arbitrary) constants satisfying:*

$$CB \exp(\sigma B) < 1 \tag{2.8}$$

*Then there exists a non-decreasing function* $Q:\Re_{+} \to \Re_{+}$ *such that for every partition* $\pi = \{\tau_i\}_{i=0}^{\infty}$ *of* $\Re_{+}$ *with* $\sup_{i \geq 0}(\tau_{i+1} - \tau_i) \leq B$ *and* $\inf_{i \geq 0}(\tau_{i+1} - \tau_i) \geq b$, *for every* $(z_0,w_0,u,v) \in \Re^{l} \times \Re^{k} \times L^{\infty}([-r,+\infty);U) \times L^{\infty}_{loc}(\Re_{+};\Re^{k})$, $(x_0,\xi_0) \in C^{0}([-r,0];\Re^{n}) \times C^{0}([-r,0];\Re^{q})$ *the unique solution of the system (1.1) with*

$$\begin{aligned} \dot{z}(t) &= F(z(t), w(t), u(t-r)) \\ \hat{x}(t) &= \Psi(z(t)) \end{aligned} \tag{2.9}$$

$$\dot{w}(t) = L_f h(\hat{x}(t), u(t-r)),\ t \in [\tau_i, \tau_{i+1}) \tag{2.10}$$

$$w(\tau_{i+1}) = h(x(\tau_{i+1} - r)) + v(\tau_{i+1}) \tag{2.11}$$

$$\begin{aligned} \dot{\xi}(t) &= F_p\left( \xi_t, u_t, \hat{x}(t), \frac{d\hat{x}}{dt}(t) \right) \\ \tilde{x}(t) &= G(\xi_t, u_t, \hat{x}(t)) \end{aligned} \tag{2.12}$$

*with initial condition* $\xi(\theta) = (\xi_0)(\theta)$, $x(\theta) = (x_0)(\theta)$, $\theta \in [-r,0]$, $(z(0),w(0)) = (z_0,w_0)$ *corresponding to inputs* $(u,v) \in L^{\infty}([-r,+\infty);U) \times L^{\infty}_{loc}(\Re_{+};\Re^{k})$ *is defined for all* $t \geq 0$ *and satisfies the estimate:*

$$\left| \tilde{x}(t) - x(t) \right| \leq e^{-\sigma t} Q\left( \|x_0\| + \|\xi_0\| + \|u\| + |z_0| + |w_0| + \sup_{0 \leq s \leq t}\left( |v(s)| \right) \right) + \frac{\gamma P \exp(\sigma B)}{1 - CB\exp(\sigma B)} \sup_{0 \leq s \leq t}\left( e^{-\sigma(t-s)} |v(s)| \right),\ \forall t \geq 0 \tag{2.13}$$

**Remark 2.4:** It should be clear that the input $v \in L^{\infty}_{loc}\left(\Re_{+};\Re^{k}\right)$ is introduced in order to describe the effect of measurement errors. Next, the structure of the observer is described:

- the continuous signal $w(t)$ attempts to approximate the continuous output signal $y(t)=h(x(t-r))$ for which only sampled measurements are available. The signal $w(t)$ is updated in an impulsive way when a new measurement becomes available (at the sampling times $\pi=\{\tau_i\}_{i=0}^{\infty}$).

-the robust global exponential observer is used with $w(t)$ as input. The observer is used in order to provide an estimate $\hat{x}(t)$ of $x(t-r)$.

-the signal $\hat{x}(t)$ is used by the robust global exponential $r-$predictor for (1.1). The observer provides an estimate $\widetilde{x}(t)$ of $x(t)$.

**Proof:** Let $0<b\le B$ be constants such that (2.8) holds and let $(z_0,w_0,u,v)\in \Re^{l}\times\Re^{k}\times L^{\infty}([-r,+\infty);U)\times L^{\infty}_{loc}\left(\Re_{+};\Re^{k}\right)$, $(x_0,\xi_0)\in C^{0}([-r,0];\Re^{n})\times C^{0}([-r,0];\Re^{q})$ be arbitrary. Let $\pi=\{\tau_i\}_{i=0}^{\infty}$ be an arbitrary partition of $\Re_{+}$ with $\sup_{i\ge 0}\left(\tau_{i+1}-\tau_i\right)\le B$ and $\inf_{i\ge 0}\left(\tau_{i+1}-\tau_i\right)\ge b$.

The solution of system (1.1) with (2.9)-(2.12) with initial condition $\xi(\theta)=(\xi_0)(\theta)$, $x(\theta)=(x_0)(\theta)$, $\theta\in[-r,0]$, $(z(0),w(0))=(z_0,w_0)$ corresponding to inputs $(u,v)\in L^{\infty}([-r,+\infty);U)\times L^{\infty}_{loc}\left(\Re_{+};\Re^{k}\right)$ exists for all $t\ge 0$. Indeed, for every integer $i\ge 0$ the solution $(x(t),z(t),w(t))$ of (1.1), (2.9)-(2.11) exists on $[\tau_i,\tau_{i+1}]$ by virtue of Hypothesis (H1). Therefore, the solution $(x(t),z(t),w(t))$ of (1.1), (2.9)-(2.11) exists for all $t\ge 0$. Hypothesis (H2) guarantees that the output signal $\hat{x}(t)=\Psi(z(t))$ is a function of class $A(\Re_{+};\Re^{n})$. Therefore, Definition 2.2 guarantees that the solution $\xi(t)$ of (2.12) exists for all $t\ge 0$ and satisfies the estimate:

$$\left|\widetilde{x}(t)-x(t)\right|\le e^{-\sigma t}a\left(\|x_0\|+\|\xi_0\|+\|u\|+|\hat{x}(0)|+\sup_{0\le s\le r}\left|\frac{d\hat{x}}{ds}(s)\right|\right)+P\sup_{0\le s\le t}\left(e^{-\sigma(t-s)}\left|\hat{x}(s)-x(s-r)\right|\right),\ \forall t\ge 0 \tag{2.14}$$

Definition (2.1) guarantees that the following estimate holds for all $t\ge r$:

$$\left|\hat{x}(t)-x(t-r)\right|\le e^{-\sigma(t-r)}M\left(\|x_0\|+|z(r)|+\|u\|\right)+\gamma\sup_{r\le s\le t}\left(e^{-\sigma(t-s)}\left|w(s)-h(x(s-r))\right|\right) \tag{2.15}$$

and hypothesis (H3) guarantees that the following estimate holds for all $t\ge r+T(x(0),z(r))$:

$$\left|L_f h(\hat{x}(t),u(t-r))-L_f h(x(t-r),u(t-r))\right|\le e^{-\sigma(t-r)}N\left(\|x_0\|+|z(r)|+\|u\|\right)+C\sup_{r\le s\le t}\left(e^{-\sigma(t-s)}\left|w(s)-h(x(s-r))\right|\right) \tag{2.16}$$

Moreover, for all $t\in[\tau_i,\tau_{i+1})$ with $\tau_i\ge r$, we obtain from (2.10) and (2.11):

$$\left|w(t)-h(x(t-r))\right|\le\left|v(\tau_i)\right|+\left(t-\tau_i\right)\sup_{\tau_i\le s\le t}\left|L_f h(\hat{x}(s),u(s-r))-L_f h(x(s-r),u(s-r))\right|$$

Consequently, since $t-\tau_i\le B$, we get from the above inequality for all $t\ge r+B$:

$$\sup_{r+B\le s\le t}\left(\left|w(s)-h(x(s-r))\right|e^{\sigma s}\right)\le e^{\sigma B}\sup_{r\le s\le t}\left(|v(s)|e^{\sigma s}\right)+Be^{\sigma B}\sup_{r\le s\le t}\left(e^{\sigma s}\left|L_f h(\hat{x}(s),u(s-r))-L_f h(x(s-r),u(s-r))\right|\right) \tag{2.17}$$

For all $t\ge r+\max(T(x(0),z(r)),B)$ it holds that $\sup_{r\le s\le t}\left(e^{\sigma s}\left|w(s)-h(x(s-r))\right|\right)=\sup_{r\le s\le r+B}\left(e^{\sigma s}\left|w(s)-h(x(s-r))\right|\right)$ or $\sup_{r\le s\le t}\left(e^{\sigma s}\left|w(s)-h(x(s-r))\right|\right)=\sup_{r+B\le s\le t}\left(e^{\sigma s}\left|w(s)-h(x(s-r))\right|\right)$. Therefore, we get from (2.16) and (2.17) for all $t\ge r+\max(T(x(0),z(r)),B)$:

$$\sup_{r+\max(B,T(x(0),z(r)))\leq s\leq t}\left(e^{\sigma s}\left|L_f h(\hat{x}(s),u(s-r))-L_f h(x(s-r),u(s-r))\right|\right)\leq$$
$$e^{\sigma r}N\left(\|x_0\|+|z(r)|+\|u\|\right)+Ce^{\sigma B}\sup_{r\leq s\leq t}\left(e^{\sigma s}|v(s)|\right) \tag{2.18}$$
$$+CBe^{\sigma B}\sup_{r\leq s\leq t}\left(e^{\sigma s}\left|L_f h(\hat{x}(s),u(s-r))-L_f h(x(s-r),u(s-r))\right|\right)$$

or

$$\sup_{r+\max(B,T(x(0),z(r)))\leq s\leq t}\left(e^{\sigma s}\left|L_f h(\hat{x}(s),u(s-r))-L_f h(x(s-r),u(s-r))\right|\right)\leq$$
$$e^{\sigma r}N\left(\|x_0\|+|z(r)|+\|u\|\right)+C\sup_{r\leq s\leq r+B}\left(e^{\sigma s}|w(s)-h(x(s-r))|\right) \tag{2.19}$$

Define $R:=r+\max(T(x(0),z(r)),B)$. Notice that (2.8), (2.18) and (2.19) in conjunction with the fact that $\sup_{r\leq s\leq t}\left(e^{\sigma s}\left|L_f h(\hat{x}(s),u(s-r))-L_f h(x(s-r),u(s-r))\right|\right)=\sup_{r\leq s\leq R}\left(e^{\sigma s}\left|L_f h(\hat{x}(s),u(s-r))-L_f h(x(s-r),u(s-r))\right|\right)$ or $\sup_{r\leq s\leq t}\left(e^{\sigma s}\left|L_f h(\hat{x}(s),u(s-r))-L_f h(x(s-r),u(s-r))\right|\right)=\sup_{R\leq s\leq t}\left(e^{\sigma s}\left|L_f h(\hat{x}(s),u(s-r))-L_f h(x(s-r),u(s-r))\right|\right)$ give for all $t\geq r$:

$$\sup_{r\leq s\leq t}\left(e^{\sigma s}\left|L_f h(\hat{x}(s),u(s-r))-L_f h(x(s-r),u(s-r))\right|\right)\leq$$
$$\frac{e^{\sigma r}}{1-CBe^{\sigma B}}N\left(\|x_0\|+|z(r)|+\|u\|\right)+\frac{Ce^{\sigma B}}{1-CBe^{\sigma B}}\sup_{r\leq s\leq t}\left(e^{\sigma s}|v(s)|\right)$$
$$+\sup_{r\leq s\leq R}\left(e^{\sigma s}\left|L_f h(\hat{x}(s),u(s-r))-L_f h(x(s-r),u(s-r))\right|\right) \tag{2.20}$$
$$+C\sup_{r\leq s\leq r+B}\left(e^{\sigma s}|w(s)-h(x(s-r))|\right)$$

Inequality (2.17) in conjunction with (2.20) implies that the following inequality holds for all $t\geq r$:

$$\sup_{r\leq s\leq t}\left(e^{\sigma s}|w(s)-h(x(s-r))|\right)\leq e^{\sigma B}\sup_{r\leq s\leq t}\left(e^{\sigma s}|v(s)|\right)$$
$$+\frac{Be^{\sigma(r+B)}}{1-CBe^{\sigma B}}N\left(\|x_0\|+|z(r)|+\|u\|\right)+\frac{CBe^{2\sigma B}}{1-CBe^{\sigma B}}\sup_{r\leq s\leq t}\left(e^{\sigma s}|v(s)|\right)$$
$$+Be^{\sigma B}\sup_{r\leq s\leq R}\left(e^{\sigma s}\left|L_f h(\hat{x}(s),u(s-r))-L_f h(x(s-r),u(s-r))\right|\right) \tag{2.21}$$
$$+(1+Be^{\sigma B}C)\sup_{r\leq s\leq r+B}\left(e^{\sigma s}|w(s)-h(x(s-r))|\right)$$

Inequality (2.15) in conjunction with (2.21) implies that the following inequality holds for all $t\geq r$:

$$\sup_{r\leq s\leq t}\left(e^{\sigma s}|\hat{x}(s)-x(s-r)|\right)\leq e^{\sigma r}M\left(\|x_0\|+|z(r)|+\|u\|\right)+\frac{\gamma e^{\sigma B}}{1-CBe^{\sigma B}}\sup_{r\leq s\leq t}\left(e^{\sigma s}|v(s)|\right)$$
$$+\gamma\frac{Be^{\sigma(r+B)}}{1-CBe^{\sigma B}}N\left(\|x_0\|+|z(r)|+\|u\|\right)+\gamma(1+Be^{\sigma B}C)\sup_{r\leq s\leq r+B}\left(e^{\sigma s}|w(s)-h(x(s-r))|\right) \tag{2.22}$$
$$+\gamma Be^{\sigma B}\sup_{r\leq s\leq R}\left(e^{\sigma s}\left|L_f h(\hat{x}(s),u(s-r))-L_f h(x(s-r),u(s-r))\right|\right)$$

Inequality (2.14) in conjunction with (2.22) implies that the following inequality holds for all $t\geq 0$:

$$\sup_{0\le s\le t}\left(e^{\sigma s}\left|\tilde{x}(s)-x(s)\right|\right)\le a\left(\|x_0\|+\|\xi_0\|+|\hat{x}(0)|+\|u\|+\sup_{0\le s\le r}\left|\frac{d\hat{x}}{ds}(s)\right|\right)$$
$$+Pe^{\sigma r}M\left(\|x_0\|+|z(r)|+\|u\|\right)+P\gamma\frac{Be^{\sigma(r+B)}}{1-CBe^{\sigma B}}N\left(\|x_0\|+|z(r)|+\|u\|\right)$$
$$+\frac{P\gamma e^{\sigma B}}{1-CBe^{\sigma B}}\sup_{r\le s\le \max(r,t)}\left(e^{\sigma s}|v(s)|\right)+P\gamma(1+Be^{\sigma B}C)\sup_{r\le s\le r+B}\left(e^{\sigma s}|w(s)-h(x(s-r))|\right)+P\sup_{0\le s\le r}\left(e^{\sigma s}|\hat{x}(s)-x(s-r)|\right)$$
$$+P\gamma Be^{\sigma B}\sup_{r\le s\le R}\left(e^{\sigma s}\left|L_f h(\hat{x}(s),u(s-r))-L_f h(x(s-r),u(s-r))\right|\right)\tag{2.23}$$

Suppose that $v\in L^\infty\left(\Re_+;\Re^k\right)$. Using Lemma 2.2 in [4], hypothesis (H1) and the fact that the mappings $\Psi(z)$, $D\Psi(z)F(z,w,u)$ are continuous, we can guarantee the existence of a non-decreasing function $c:\Re_+\to\Re_+$ such that for every $t\in[\tau_i,\tau_{i+1})$ the following estimate holds:

$$|z(t)|+|w(t)|+|\Psi(z(t))|+|D\Psi(z(t))F(z(t),w(t),u(t-r))|\le c\left(|z(\tau_i)|+|w(\tau_i)|+\|u\|\right)\tag{2.24}$$

It follows from (2.24) and the fact that $\inf_{i\ge 0}\left(\tau_{i+1}-\tau_i\right)\ge b$ (which directly implies that at most $N=\left[\frac{r}{b}\right]$ points of the partition $\pi=\{\tau_i\}_{i=0}^\infty$ are in the interval $[0,r)$) that there exists a non-decreasing function $\tilde{c}:\Re_+\to\Re_+$ such that

$$\sup_{0\le s\le r}|\hat{x}(s)|+|z(r)|+\sup_{0\le s\le r}\left|\frac{d\hat{x}}{ds}(s)\right|\le\tilde{c}\left(|z_0|+|w_0|+\|u\|+\|x_0\|+\|v\|\right)\tag{2.25}$$

Since $T(x,z)$ is continuous, inequality (2.25) implies that there exists a non-decreasing function $\phi:\Re_+\to\Re_+$ such that

$$R=r+\max(B,T(x(0),z(r)))\le\phi\left(|z_0|+|w_0|+\|u\|+\|x_0\|+\|v\|\right)\tag{2.26}$$

It follows from Lemma 2.2 in [4], the fact that (1.1) is forward complete, (2.24) and the fact that $\inf_{i\ge 0}\left(\tau_{i+1}-\tau_i\right)\ge b$ (which directly implies that at most $N=\left[\frac{\phi\left(|z_0|+|w_0|+\|u\|+\|x_0\|+\|v\|\right)}{b}\right]$ points of the partition $\pi=\{\tau_i\}_{i=0}^\infty$ are in the interval $\left[0,\phi\left(|z_0|+|w_0|+\|u\|+\|x_0\|+\|v\|\right)\right)$) that there exists a non-decreasing function $\bar{c}:\Re_+\to\Re_+$ such that

$$\sup_{r\le s\le R}\left|L_f h(\hat{x}(s),u(s-r))-L_f h(x(s-r),u(s-r))\right|+\sup_{r\le s\le R}|w(s)-h(x(s-r))|\le\bar{c}\left(|z_0|+|w_0|+\|u\|+\|x_0\|+\|v\|\right)\tag{2.27}$$

Using (2.23), (2.25), (2.26), (2.27), it follows that there exists a non-decreasing function $Q:\Re_+\to\Re_+$ such that the following inequality holds for all $t\ge 0$:

$$\sup_{0\le s\le t}\left(e^{\sigma s}|\tilde{x}(s)-x(s)|\right)\le Q\left(\|x_0\|+\|\xi_0\|+|w_0|+|z_0|+\|u\|+\|v\|\right)+\frac{P\gamma e^{\sigma B}}{1-CBe^{\sigma B}}\sup_{r\le s\le\max(r,t)}\left(e^{\sigma s}|v(s)|\right)\tag{2.28}$$

Inequality (2.13) is a direct consequence of (2.28) and the causality property for system (1.1) with (2.9)-(2.12). The proof is complete. $\lhd$

## 3. Globally Lipschitz Systems

In this section we consider the construction of global exponential sampled-data observers for globally Lipschitz systems. We consider system (1.1), (1.2) and we assume that

**(H4)** *There exists a constant* $L>0$ *such that*

$$|f(x,u)-f(z,u)| \le L|x-z| \, , \forall x,z \in \Re^n \, , \; \forall u \in U \tag{3.1}$$

**(H5)** *There exists a symmetric, positive definite matrix* $P \in \Re^{n \times n}$*, a constant* $q>0$ *and matrices* $K \in \Re^{n \times k}$*,* $H \in \Re^{k \times n}$ *such that:*

$$h(x)=Hx \, , \; \forall x \in \Re^n \tag{3.2}$$

$$(z-x)'P\big(f(z,u)-f(x,u)\big)+(z-x)'PKH(z-x) \le -q|z-x|^2 \, , \; \forall x,z \in \Re^n \, , \; \forall u \in U \tag{3.3}$$

Hypotheses (H4), (H5) are automatically satisfied for triangular systems of the form:

$$\begin{aligned}
&\dot{x}_i = f_i(u,x_1,...,x_i)+x_{i+1} \quad , i=1,...,n-1 \\
&\dot{x}_n = f_n(u,x_1,...,x_n) \\
&y = x_1 \in \Re \\
&x=(x_1,...,x_n) \in \Re^n
\end{aligned}$$

where the smooth mappings $f_i : U \times \Re^i \to \Re$ ($i=1,...,n$) are globally Lipschitz with respect to $x \in \Re^n$ (see [9]).

Hypothesis (H5) guarantees that (2.1) with $l=n$, $F(z,y,u)=f(z,u)+K\big(Hz-y\big)$ and $\Psi(z)=z \in \Re^n$ is a Robust Global Exponential Observer for (1.1). Moreover, due to hypotheses (H4), (H5), the system $\dot{z}=f(z,u)+K\big(Hz-w\big)$ with inputs $u,w$ and the system $\dot{z}=f(z,u)+K\big(Hz-w\big)$, $\dot{w}=Hf(z,u)$ with input $u$, are forward complete. Consequently, hypothesis (H1) holds.

Notice that, by virtue of (3.3), for every $(x_0,z_0,u,v) \in \Re^n \times \Re^n \times L^\infty\big(\Re_+;U\big) \times L^\infty_{loc}\big(\Re_+;\Re^k\big)$ the solution $(x(t),z(t))$ of (1.1), (1.2) and (2.1) with initial condition $(x(0),z(0))=(x_0,z_0)$ corresponding to inputs $(u,v) \in L^\infty\big(\Re_+;U\big) \times L^\infty_{loc}\big(\Re_+;\Re^k\big)$ satisfies the following estimate for all $t \ge 0$:

$$|\hat{x}(t)-x(t)| \le \sqrt{\frac{|P|}{R}} e^{-\sigma t}|z_0-x_0| + \frac{\sqrt{|K'PPK|}}{q}\sqrt{\frac{|P|}{R}} \sup_{0 \le s \le t}\left(e^{-\sigma(t-s)}|v(s)|\right) \tag{3.4}$$

where $R:=\min\{x'Px : |x|=1\}$ and $\sigma := \frac{q}{2|P|}$. Therefore, by virtue of (3.1), (3.2) and (3.4), it follows that Hypothesis (H3) holds with $\sigma := \frac{q}{2|P|}$, $T(x_0,z_0) \equiv 0$, $C := L|H|\frac{\sqrt{|K'PPK|}}{q}\sqrt{\frac{|P|}{R}}$ and $N(s) := s\sqrt{\frac{|P|}{R}}$.

We will show next that Hypothesis (H2) holds as well with $A(\Re_+;\Re^n)$ being the set of the functions $z:\Re_+ \to \Re^n$ which are absolutely continuous on every bounded interval of $\Re_+$ with $\dot{z} \in L^\infty_{loc}\big(\Re_+;\Re^n\big)$. The following theorem guarantees (in a constructive way) that there exists a robust global exponential $r-$predictor for (1.1).

**Theorem 3.1:** *Consider system (1.1) under hypothesis (H4) and let* $\sigma,r>0$ *be given constants. For every constant* $\mu>\sigma$ *and for every positive integer* $p \ge 1$ *with* $L\frac{\exp\left(\sigma r p^{-1}\right)-1}{\sigma}<1$ *there exist constants* $Q_j>0$ *(* $j=1,...,5$ *), such that for every* $u \in L^\infty_{loc}\big([-r,+\infty);U\big)$*, for every absolutely continuous mapping* $z:\Re_+ \to \Re^n$ *with* $\dot{z} \in L^\infty_{loc}\big(\Re_+;\Re^n\big)$ *and for every* $x_0 \in C^0([-r,0];\Re^n)$*,* $\xi_{i,0} \in C^0([-r,0];\Re^n)$ *(* $i=1,...,p$ *), the unique solution of system (1.1) with*

$$\dot{\xi}_1(t) = f(\xi_1(t), u(t-r+\delta)) - f(\xi_1(t-\delta), u(t-r)) + \dot{z}(t) - \mu\left( \xi_1(t) - z(t) - \int_{t-\delta}^{t} f\left(\xi_1(s), u(s-r+\delta)\right) ds \right) \quad (3.5)$$

$$\begin{aligned} \dot{\xi}_i(t) &= f(\xi_i(t), u(t-r+i\delta)) - f(\xi_i(t-\delta), u(t-r+(i-1)\delta)) \\ &+ \dot{\xi}_{i-1}(t) - \mu\left( \xi_i(t) - \xi_{i-1}(t) - \int_{t-\delta}^{t} f\left(\xi_i(s), u(s-r+i\delta)\right) ds \right), \; i = 2,...,p \end{aligned} \quad (3.6)$$

*where* $\delta := p^{-1} r$ *, with initial condition* $x(\theta) = x_0(\theta)$ *,* $\xi_i(\theta) = \xi_{i,0}(\theta)$ *,* $\theta \in [-r, 0]$ *, (*$i = 1,...,p$*), exists for all* $t \geq 0$ *and satisfies the following estimate for all* $i = 1,...,p$ *and* $t \geq 0$ *:*

$$\begin{aligned} &\left| \xi_i(t) - x(t-r+i\delta) \right| \\ &\leq \beta^i \sup_{0 \leq s \leq t} \left( \exp(-\sigma(t-s)) \left| z(s) - x(s-r) \right| \right) + Q_1 \exp(-\sigma t) \sup_{-r \leq s \leq 0} \left( |x(s)| \right) \\ &+ Q_2 \exp(-\sigma t) \sup_{-r \leq s \leq r+\delta} \left( |f(0, u(s))| \right) + Q_3 \exp(-\sigma t) \sum_{k=1}^{i} \sup_{-\delta \leq s \leq 0} \left( |\xi_k(s)| \right) \\ &+ Q_4 \exp(-\sigma t) |z(0)| + Q_5 \exp(-\sigma t) \sup_{0 \leq s \leq r} \left( |\dot{z}(s)| \right) \end{aligned} \quad (3.7)$$

*where* $\beta := \dfrac{\sigma}{\sigma - L(\exp(\sigma\delta) - 1)}$.

**Remark 3.2:** It is clear that inequality (3.7) guarantees that the system (3.5), (3.6) with output $y(t) = \xi_p(t)$ is a robust global exponential $r$ – predictor for (1.1).

The proof of Theorem 3.1 is an inductive application of the following technical lemma, which is proved at the Appendix.

**Lemma 3.3:** *Consider system (1.1) under hypothesis (H4) and let* $\sigma, r > 0$ *be constants. For every constant* $\mu > \sigma$ *and for every* $\delta \in [0, r]$ *with* $L \dfrac{\exp(\sigma\delta) - 1}{\sigma} < 1$ *there exist constants* $Q_i > 0$ *(*$i = 1,...,5$*), such that for every* $u \in L_{loc}^{\infty}([-r, +\infty); U)$*, for every absolutely continuous mapping* $z : \Re_+ \to \Re^n$ *with* $\dot{z} \in L_{loc}^{\infty}(\Re_+; \Re^n)$ *and for every* $x_0 \in C^0([-r, 0]; \Re^n)$*,* $\xi_0 \in C^0([-r, 0]; \Re^n)$*, the unique solution of system (1.1) with*

$$\dot{\xi}(t) = f(\xi(t), u(t-r+\delta)) - f(\xi(t-\delta), u(t-r)) + \dot{z}(t) - \mu\left( \xi(t) - z(t) - \int_{t-\delta}^{t} f\left(\xi(s), u(s-r+\delta)\right) ds \right) \quad (3.8)$$

*with initial condition* $x(\theta) = x_0(\theta)$*,* $\xi(\theta) = \xi_0(\theta)$*,* $\theta \in [-r, 0]$*, exists for all* $t \geq 0$ *and satisfies the following estimate for all* $t \geq 0$ *:*

$$\begin{aligned} &\left| \xi(t) - x(t-r+\delta) \right| \\ &\leq \beta \sup_{0 \leq s \leq t} \left( \exp(-\sigma(t-s)) \left| z(s) - x(s-r) \right| \right) + Q_1 \exp(-\sigma t) \sup_{-r \leq s \leq 0} \left( |x(s)| \right) \\ &+ Q_2 \exp(-\sigma t) \sup_{-r \leq s \leq \delta} \left( |f(0, u(s))| \right) + Q_3 \exp(-\sigma t) \sup_{-\delta \leq s \leq 0} \left( |\xi(s)| \right) \\ &+ Q_4 \exp(-\sigma t) |z(0)| + Q_5 \exp(-\sigma t) \sup_{0 \leq s \leq r} \left( |\dot{z}(s)| \right) \end{aligned} \quad (3.9)$$

*where* $\beta := \dfrac{\sigma}{\sigma - L(\exp(\sigma\delta) - 1)}$.

Therefore, all the above enable us to design a robust global sampled-data exponential observer with delayed measurements.

**Theorem 3.4:** *Consider system (1.1) under hypotheses (H4-5) and let* $r>0$ *be constant. Then the following system*

$$\dot{z}(t) = f(z(t),u(t-r)) + K(Hz(t)-w(t)) \tag{3.10}$$

$$\dot{w}(t) = Hf(z(t),u(t-r))\,,\ t\in[\tau_i,\tau_{i+1}) \tag{3.11}$$

$$w(\tau_{i+1}) = Hx(\tau_{i+1}-r) + v(\tau_{i+1}) \tag{3.12}$$

*with (3.5), (3.6),* $\mu>0$, $\delta := p^{-1}r$ *and* $L\delta<1$, *is a robust global exponential sampled-data observer for (1.1), provided that the upper diameter* $B$ *of the sampling partition satisfies* $L|H|\frac{\sqrt{|K'PPK|}}{q}\sqrt{\frac{|P|}{R}}B<1$, *i.e., for every positive integer* $p>0$ *with* $Lr<p$, *for every* $\mu>0$, $0<b\le B$ *with* $L|H|\frac{\sqrt{|K'PPK|}}{q}\sqrt{\frac{|P|}{R}}B<1$, *there exist a non-decreasing function* $Q:\Re_+\to\Re_+$ *and constants* $\sigma,\Gamma>0$ *such that for every partition* $\pi=\{\tau_i\}_{i=0}^{\infty}$ *of* $\Re_+$ *with* $\sup_{i\ge0}(\tau_{i+1}-\tau_i)\le B$ *and* $\inf_{i\ge0}(\tau_{i+1}-\tau_i)\ge b$, *for every* $(z_0,w_0,u,v)\in\Re^n\times\Re^k\times L^\infty([-r,+\infty);U)\times L^\infty_{loc}(\Re_+;\Re^k)$, $x_0\in C^0([-r,0];\Re^n)$, $\xi_{i,0}\in C^0([-r,0];\Re^n)$ $(i=1,...,p)$ *the unique solution of the system (1.1) with (3.5), (3.6), (3.10), (3.11), (3.12) with* $\delta := p^{-1}r$, *initial condition* $x(\theta)=x_0(\theta)$, $\xi_i(\theta)=\xi_{i,0}(\theta)$, $\theta\in[-r,0]$, $(i=1,...,p)$, $(z(0),w(0))=(z_0,w_0)$ *corresponding to inputs* $(u,v)\in L^\infty([-r,+\infty);U)\times L^\infty_{loc}(\Re_+;\Re^k)$ *is defined for all* $t\ge0$ *and satisfies the estimate:*

$$\left|\xi_p(t)-x(t)\right| \le e^{-\sigma t}Q\left(\|x_0\|+\sum_{i=1}^{p}\|\xi_{i,0}\|+\|u\|+|z_0|+|w_0|+\sup_{0\le s\le t}\left(|v(s)|\right)\right)+\Gamma\sup_{0\le s\le t}\left(e^{-\sigma(t-s)}|v(s)|\right),\ \forall t\ge0 \tag{3.13}$$

**Proof:** If the above inequalities hold then (by continuity) there exists $\sigma>0$ (sufficiently small) such that $L|H|\frac{\sqrt{|K'PPK|}}{p}\sqrt{\frac{|P|}{R}}B\exp(\sigma B)<1$, $L\frac{\exp(\sigma\delta)-1}{\sigma}<1$ and such that all Hypotheses (H1-3) hold. The rest is a direct consequence of Theorem 2.3. $\lhd$

A direct application to the Linear Time-Invariant case $\dot{x}=Fx+Gu$ (where $F\in\Re^{n\times n}$, $G\in\Re^{n\times m}$) guarantees that the linear system:

$$\dot{\xi}_0(t) = (F+KH)\xi_0(t) + Gu(t-r) - Kw(t) \tag{3.14}$$

$$\dot{w}(t) = HF\xi_0(t) + HGu(t-r)\,,\ t\in[\tau_i,\tau_{i+1}) \tag{3.15}$$

$$w(\tau_{i+1}) = Hx(\tau_{i+1}-r) + v(\tau_{i+1}) \tag{3.16}$$

$$\begin{aligned}\dot{\xi}_i(t) &= F\left(\xi_i(t)-\xi_i(t-\delta)+\mu\int_{t-\delta}^{t}\xi_i(s)ds\right)+\dot{\xi}_{i-1}(t)-\mu\left(\xi_i(t)-\xi_{i-1}(t)\right)\\ &+G\left(u(t-r+i\delta)-u(t-r+(i-1)\delta)+\mu\int_{t-\delta}^{t}u(s-r+i\delta)ds\right)\end{aligned}\ ,\ i=1,...,p \tag{3.17}$$

with $\mu>0$, $\delta := p^{-1}r$ and $|F|r<p$, is a robust global exponential sampled-data observer for $\dot{x}=Fx+Gu$, provided that the upper diameter $B$ of the sampling partition satisfies $|F||H|\frac{\sqrt{|K'PPK|}}{q}\sqrt{\frac{|P|}{R}}B<1$.

## 4. Systems with a Globally Asymptotically Stable Set

In this section we consider the construction of global exponential sampled-data observers for nonlinear systems with a globally asymptotically stable set. More specifically, we consider system (1.1), (1.2) and we assume that

**(H6)** *The set* $U \subset \Re^m$ *is compact and there exist a non-empty compact set* $S \subset \Re^n$, *a continuous function* $T : \Re^n \to \Re_+$ *and a smooth positive function* $\psi : \Re^n \to (0,+\infty)$ *such that for every* $u \in L^\infty(\Re_+;U)$ *and for every initial condition* $x(0) \in \Re^n$ *the solution of (1.1) satisfies:*

$$x(t) \in S \,, \forall t \geq T(x(0)) \tag{4.1}$$

$$|x(t)| \leq \psi(x(0)) \,, \ \forall t \geq 0 \tag{4.2}$$

For systems satisfying hypothesis (H6), a general procedure for the design of robust global exponential observers of the form (2.1) with $l = n$ was proposed in [16]. More specifically, the proof of Theorem 2.2 in [16] shows that the observer satisfies the following hypothesis:

**(H7)** $h(\Re^n) = \Re$ *and for every* $(u,w) \in L^\infty([-r,+\infty);U) \times L^\infty_{loc}([-r,+\infty);\Re^k)$ *and for every initial condition* $z(0) \in \Re^n$ *the solution of (2.6) satisfies:*

$$z(t) \in S \,, \forall t \geq T(z(0)) \tag{4.3}$$

$$|z(t)| \leq \psi(z(0)) \,, \ \forall t \geq 0 \tag{4.4}$$

Hypothesis (H7) guarantees that hypothesis (H1) holds. Moreover, hypotheses (H6), (H7) guarantee that there exists a constant $G \geq 0$ such that for every $(x_0, z_0, u, v) \in \Re^n \times \Re^l \times L^\infty(\Re_+;U) \times L^\infty_{loc}(\Re_+;\Re^k)$ the solution $(x(t), z(t))$ of (1.1), (1.2) and (2.1) with initial condition $(x(0), z(0)) = (x_0, z_0)$ corresponding to inputs $(u,v) \in L^\infty(\Re_+;U) \times L^\infty_{loc}(\Re_+;\Re^k)$ satisfies the following estimate for all $t \geq \max(T(x_0), T(z_0))$:

$$|L_f h(\hat{x}(t), u(t)) - L_f h(x(t), u(t))| \leq G |\hat{x}(t) - x(t)| \tag{4.5}$$

Inequality (4.5) in conjunction with (2.2) shows that Hypothesis (H3) holds with $N(s) := GM(s)$, $C := G\gamma$ and $T(x_0, z_0) := \max(T(x_0), T(z_0))$.

An example of a system satisfying hypotheses (H6-7) can be found in [16] (Example 4.1).

We will show next that Hypothesis (H2) holds as well with $A(\Re_+;\Re^n)$ being the set of the functions $z : \Re_+ \to \Re^n$ which are absolutely continuous on every bounded interval of $\Re_+$ with $\dot{z} \in L^\infty_{loc}(\Re_+;\Re^n)$, $|z(t)| \leq \psi(z(s))$ for all $0 \leq s \leq t$ and $z(t) \in S$, for all $t \geq T(z(0))$. The following proposition guarantees (in a constructive way) that there exists a robust global exponential $r-$predictor for (1.1). Its proof is provided at the Appendix.

**Proposition 4.1:** *Consider system (1.1) under hypothesis (H6). Let* $\sigma, r > 0$ *be constants and let* $q : \Re \to \Re_+$ *be a continuously differentiable function with* $q(s) = 1$ *for* $s \leq 1$ *and* $sq(s) \leq K$ *for* $s \geq 1$, *where* $K \geq 1$ *is a constant. Furthermore, define:*

$$p(s) := \max\left\{ |f(\xi,u)| : u \in U \,, |\xi| \leq K \max\{\psi(z) : |z| \leq s\} \right\} \tag{4.6}$$

$$a := \max\left\{ |z| : z \in S \right\} \tag{4.7}$$

$$\tilde{S} := \left\{ \xi \in \Re^n : |\xi| \leq 1 + a + r\, p(a) \right\} \tag{4.8}$$

*Let $G_1, G_2 > 0$ be constants satisfying*

$$\left| f\left( q\left(\frac{|\xi|}{\psi(z)}\right)\xi, u \right) - f\left( q\left(\frac{|x|}{\psi(y)}\right)x, u \right) \right| \le G_1|\xi - x| + G_2|z - y| \qquad (4.9)$$
$$\forall u \in U,\, y, x, z \in S,\, \xi \in \widetilde{S}$$

*Then for every constant $\mu \ge \sigma$, for every $\delta \in [0, r]$ with $G_1 \frac{\exp(\sigma\delta) - 1}{\sigma} < 1$, there exist a non-decreasing function $M : \Re_+ \to \Re_+$, such that for every $u \in L^\infty_{loc}([-r, +\infty); U)$, for every absolutely continuous mapping $z : \Re_+ \to \Re^n$ with $\dot{z} \in L^\infty_{loc}(\Re_+; \Re^n)$, $|z(t)| \le \psi(z(s))$ for all $0 \le s \le t$ and $z(t) \in S$, for all $t \ge T(z(0))$ and for every $x_0 \in C^0([-r, 0]; \Re^n)$, $\xi_0 \in C^0([-r, 0]; \Re^n)$, the unique solution of system (1.1) with*

$$\dot{\xi}(t) = \dot{z}(t) + \frac{d}{dt}\int_{t-\delta}^{t} f\left( q\left(\frac{|\xi(s)|}{\psi(z(t))}\right)\xi(s), u(s - r + \delta) \right) ds - \mu\left( \xi(t) - z(t) - \int_{t-\delta}^{t} f\left( q\left(\frac{|\xi(s)|}{\psi(z(t))}\right)\xi(s), u(s - r + \delta) \right) ds \right) \qquad (4.10)$$

*with initial condition $x(\theta) = x_0(\theta)$, $\xi(\theta) = \xi_0(\theta)$, $\theta \in [-r, 0]$, exists for all $t \ge 0$ and satisfies the following estimate for all $t \ge 0$:*

$$|\xi(t) - x(t - r + \delta)| \le \exp(-\sigma t) M\left( \|x_0\| + |z(0)| + \|\xi_0\| \right) + \beta \sup_{0 \le s \le t}\left( \exp(-\sigma(t - s)) |z(s) - x(s - r)| \right) \qquad (4.11)$$

*where $\beta := \dfrac{\sigma + G_2(\exp(\sigma\delta) - 1)}{\sigma - G_1(\exp(\sigma\delta) - 1)}$.*

**Remark 4.2:** An example of a function $q : \Re \to \Re_+$ that satisfies the requirements of Proposition 4.1 is the function $q(s) := 2s^{-1} - s^{-2}$ for $s > 1$ and $q(s) = 1$ for $s \le 1$.

We are now in a position to state and prove the main result of this section.

**Theorem 4.3:** *Consider system (1.1) under hypotheses (H6-7) and let $r > 0$ be a constant. If $G_1 r < 1$, where $G_1 > 0$ is the constant involved in (4.9) then for every $\mu > 0$, $0 < b \le B$ with $G\gamma B < 1$, where $G > 0$ is the constant involved in (4.5) and $\gamma$ is the constant involved in (2.2), there exist a non-decreasing function $Q : \Re_+ \to \Re_+$ and constants $\sigma, \Gamma > 0$ such that for every partition $\pi = \{\tau_i\}_{i=0}^{\infty}$ of $\Re_+$ with $\sup_{i \ge 0}(\tau_{i+1} - \tau_i) \le B$ and $\inf_{i \ge 0}(\tau_{i+1} - \tau_i) \ge b$, for every $(z_0, w_0, u, v) \in \Re^n \times \Re^k \times L^\infty([-r, +\infty); U) \times L^\infty_{loc}(\Re_+; \Re^k)$, $(x_0, \xi_0) \in C^0([-r, 0]; \Re^n) \times C^0([-r, 0]; \Re^n)$ the unique solution of the system (1.1), (2.9), (2.10), (2.11), (4.10) with $\delta = r$, initial condition $x(\theta) = x_0(\theta)$, $\xi(\theta) = \xi_0(\theta)$, $\theta \in [-r, 0]$, $(z(0), w(0)) = (z_0, w_0)$ corresponding to inputs $(u, v) \in L^\infty([-r, +\infty); U) \times L^\infty_{loc}(\Re_+; \Re^k)$ is defined for all $t \ge 0$ and satisfies the estimate:*

$$|\xi(t) - x(t)| \le e^{-\sigma t} Q\left( \|x_0\| + \|\xi_0\| + \|u\| + |z_0| + |w_0| + \sup_{0 \le s \le t}(|v(s)|) \right) + \Gamma \sup_{0 \le s \le t}\left( e^{-\sigma(t - s)}|v(s)| \right),\ \forall t \ge 0 \qquad (4.12)$$

**Proof:** If $G_1 r < 1$ then (by continuity) there exists sufficiently small $\sigma > 0$ such that $G_1 \frac{\exp(\sigma\delta) - 1}{\sigma} < 1$ and hypothesis (H2) holds with state predictor defined by (4.10) with $\delta = r$, $\mu > 0$. Using hypotheses (H6-7), we deduce that (H1) and (H3) also hold as explained above. If the upper diameter $B$ of the sampling partition satisfies $G\gamma B < 1$ then, by continuity there exists sufficiently small $\sigma > 0$ such that $G\gamma B \exp(\sigma B) < 1$. The rest is a direct consequence of Theorem 2.3. $\triangleleft$

## 5. Concluding remarks

This present work provides new results concerning:

1) the robust global exponential sampled-data observer design for wide classes of nonlinear systems with both sampled and delayed measurements (Theorems 2.3, 3.4 and 4.3), and

2) the robust global exponential state predictor design for wide classes of nonlinear systems (Theorem 3.1 and Proposition 4.1).

The global exponential state predictors are constructed by means of small-gain arguments and additional conditions on the prediction time horizon. It is also shown that, if a special structure cascade is used with a sufficient number of predictors, then an exponential state predictor for an arbitrary prediction time horizon can be constructed for the special class of globally Lipschitz systems (Theorem 3.1).

The global exponential sampled-data observer design is accomplished by using a small gain approach and sufficient conditions are provided, which involve both the delay and the sampling period. The structure of the proposed observer can be described as follows:
- a hybrid sampled-data observer is first used in order to utilize the sampled and delayed measurements and provide an estimate of the delayed state vector,
-the estimate of the delayed state vector is used by the robust global exponential predictor. The predictor provides an estimate of the current value of the state vector.

The proposed robust global exponential sampled-data observer is robust with respect to measurement errors and perturbations of the sampling schedule.

The obtained results can be used in a straightforward way for the stabilization of nonlinear systems with input delays. The predictor-based feedback, which was recently proposed in [17,18,19,20], can be used in conjunction with the proposed robust global exponential sampled-data observers in order to solve the output feedback stabilization problem for nonlinear systems with input delays. This will be the topic of future research.

# Appendix

**Proof of Lemma 3.3:** First notice that by virtue of (3.1), the right hand side of (3.8) satisfies the inequality:

$$\left|\dot{\xi}(t)\right| \leq \left|\dot{z}(t)+\mu z(t)\right| + \left(\mu L\delta+\mu+2L\right)\sup_{t-\delta\leq s\leq t}\left(\left|\xi(s)\right|\right) + \left(\mu\delta+2\right)\sup_{t-\delta\leq s\leq t}\left(\left|f\left(0,u(s-r+\delta)\right)\right|\right) \tag{A.1}$$

for almost all $t \geq 0$ for which the solution of (3.8) exists. Notice that the continuous functional $V(t)=\sup_{t-\delta\leq s\leq t}\left(\left|\xi(s)\right|\right)$ satisfies:

$$\limsup_{h\to 0^+}\frac{V(t+h)-V(t)}{h} \leq \left(\mu L\delta+\mu+2L\right)V(t) + \sup_{0\leq s\leq t}\left|\dot{z}(s)+\mu z(s)\right| + \left(\mu\delta+2\right)\sup_{-r\leq s\leq t-r+\delta}\left|f(0,u(s))\right| \tag{A.2}$$

Using the comparison lemma (Lemma 2.12, page 77 in [15]), we obtain:

$$\sup_{t-\delta\leq s\leq t}\left|\xi(s)\right| \leq \exp\left(\left(\mu L\delta+\mu+2L\right)t\right)\left(\sup_{-\delta\leq s\leq 0}\left|\xi(s)\right| + \frac{\sup_{0\leq s\leq t}\left|\dot{z}(s)+\mu z(s)\right| + \left(\mu\delta+2\right)\sup_{-r\leq s\leq t-r+\delta}\left|f(0,u(s))\right|}{\mu L\delta+\mu+2L}\right) \tag{A.3}$$

The above inequality shows that the solution of (3.8) exists for all $t \geq 0$. Next, notice that the solution of (3.8) satisfies the following equation for all $t \geq 0$:

$$\xi(t)-z(t)-\int_{t-\delta}^{t} f\left(\xi(s),u(s-r+\delta)\right)ds = \exp\left(-\mu t\right)\left(\xi(0)-z(0)-\int_{-\delta}^{0} f\left(\xi(s),u(s-r+\delta)\right)ds\right) \tag{A.4}$$

Moreover, the following equation holds for all $t \geq r$ :

$$x(t-r+\delta) = x(t-r) + \int_{t-\delta}^{t} f\big(x(s-r+\delta), u(s-r+\delta)\big)ds \tag{A.5}$$

Using (3.1), (A.4) and (A.5) we get for all $t \geq r$ :

$$\begin{aligned} |\xi(t) - x(t-r+\delta)| &\leq |z(t) - x(t-r)| + L\int_{t-\delta}^{t} |\xi(s) - x(s-r+\delta)|ds \\ &+ \exp(-\mu t)\left|\xi(0) - z(0) - \int_{-\delta}^{0} f(\xi(s), u(s-r+\delta))ds\right| \end{aligned} \tag{A.6}$$

Since $\sigma \in (0, \mu)$ we obtain from (A.6) for all $t \geq r$ :

$$\begin{aligned} \exp(\sigma t)|\xi(t) - x(t-r+\delta)| &\leq \exp(\sigma t)|z(t) - x(t-r)| \\ &+ L\frac{\exp(\sigma\delta)-1}{\sigma}\sup_{t-\delta\leq s\leq t}\big(\exp(\sigma s)|\xi(s) - x(s-r+\delta)|\big) \\ &+ \left|\xi(0) - z(0) - \int_{-\delta}^{0} f(\xi(s), u(s-r+\delta))ds\right| \end{aligned}$$

which directly implies for all $t \geq r$ :

$$\begin{aligned} \sup_{r\leq s\leq t}\big(\exp(\sigma s)|\xi(s) - x(s-r+\delta)|\big) &\leq \sup_{r\leq s\leq t}\big(\exp(\sigma s)|z(s) - x(s-r)|\big) \\ &+ L\frac{\exp(\sigma\delta)-1}{\sigma}\sup_{r-\delta\leq s\leq t}\big(\exp(\sigma s)|\xi(s) - x(s-r+\delta)|\big) \\ &+ \left|\xi(0) - z(0) - \int_{-\delta}^{0} f(\xi(s), u(s-r+\delta))ds\right| \end{aligned} \tag{A.7}$$

By distinguishing the cases $\sup_{r\leq s\leq t}\big(\exp(\sigma s)|\xi(s) - x(s-r+\delta)|\big) = \sup_{r-\delta\leq s\leq t}\big(\exp(\sigma s)|\xi(s) - x(s-r+\delta)|\big)$ and $\sup_{r-\delta\leq s\leq r}\big(\exp(\sigma s)|\xi(s) - x(s-r+\delta)|\big) = \sup_{r-\delta\leq s\leq t}\big(\exp(\sigma s)|\xi(s) - x(s-r+\delta)|\big)$, and using the inequality $L\frac{\exp(\sigma\delta)-1}{\sigma} < 1$ we obtain for all $t \geq r$ :

$$\begin{aligned} &\sup_{r\leq s\leq t}\big(\exp(\sigma s)|\xi(s) - x(s-r+\delta)|\big) \\ &\leq \frac{\sigma}{\sigma - L(\exp(\sigma\delta)-1)}\sup_{r\leq s\leq t}\big(\exp(\sigma s)|z(s) - x(s-r)|\big) \\ &+ \sup_{r-\delta\leq s\leq r}\big(\exp(\sigma s)|\xi(s) - x(s-r+\delta)|\big) + \frac{\sigma}{\sigma - L(\exp(\sigma\delta)-1)}\left|\xi(0) - z(0) - \int_{-\delta}^{0} f(\xi(s), u(s-r+\delta))ds\right| \end{aligned} \tag{A.8}$$

Using (2.1), evaluating of the Dini derivative of the continuous functional $V(t) = \sup_{t-r\leq s\leq t}\big(|x(s)|\big)$ and using the comparison lemma (Lemma 2.12, page 77 in [15]), we obtain for all $t \geq 0$ :

$$\sup_{t-r\leq s\leq t}|x(s)| \leq \exp(Lt)\left(\sup_{-r\leq s\leq 0}|x(s)| + L^{-1}\sup_{0\leq s\leq t}|f(0, u(s))|\right) \tag{A.9}$$

Combining (A.3), (A.8), (A.9) and the fact that

$$\sup_{0\le s\le r}|z(s)| \le |z(0)| + r \sup_{0\le s\le r}|\dot{z}(s)|$$

we are in a position to conclude that there exist constants $Q_i > 0$ ($i = 1,...,5$) such that (3.9) holds for all $t \ge 0$. The proof is complete. $\triangleleft$

**Proof of Proposition 4.1:** We notice that the solution of (4.10) satisfies

$$\begin{aligned} \xi(t) &= z(t) + \int_{t-\delta}^{t} f\left( q\left( \frac{|\xi(s)|}{\psi(z(t))} \right)\xi(s), u(s-r+\delta) \right) ds \\ &+ \exp(-\mu t)\left( \xi(0) - z(0) - \int_{-\delta}^{0} f\left( q\left( \frac{|\xi(s)|}{\psi(z(0))} \right)\xi(s), u(s-r+\delta) \right) ds \right) \end{aligned} \tag{A.10}$$

for all $t \ge 0$ for which the solution of (4.10) exists. Using (4.6) and (A.10) and the fact that $q\left( \frac{|\xi(s)|}{\psi(z(t))} \right)|\xi(s)| \le K\psi(z(t))$, we obtain:

$$|\xi(t)| \le |z(t)| + \delta\, p(|z(t)|) + \exp(-\mu t)\left(|\xi(0)| + |z(0)| + \delta\, p(|z(0)|)\right) \tag{A.11}$$

for all $t \ge 0$ for which the solution of (4.10) exists. A standard contradiction argument in conjunction with (A.11) shows that the solution of (4.10) exists for all $t \ge 0$. Define:

$$\tilde{T} := \max\left\{ T(x(0)) + r + \delta, T(z(0)) + \delta, \mu^{-1}\ln\left(1 + |\xi(0)| + |z(0)| + \delta\, p(|z(0)|)\right) \right\} + \delta \tag{A.12}$$

Using the facts that $z(t) \in S$ for all $t \ge T(z(0))$ and $\delta \in [0,r]$, in conjunction with (4.1), (A.11), (4.7), (4.8) and definition (A.12), we conclude that

$$\xi(t-\delta) \in \tilde{S}, \text{ for all } t \ge \tilde{T} \tag{A.13}$$

$$z(t-\delta) \in S, \text{ for all } t \ge \tilde{T} \tag{A.14}$$

$$x(t-\delta-r) \in S, \text{ for all } t \ge \tilde{T} \tag{A.15}$$

By virtue of (4.2) and the semigroup property for the solutions of (1.1), we conclude that the following equation holds for all $t \ge r$:

$$x(t-r+\delta) = x(t-r) + \int_{t-r}^{t-r+\delta} f\left( q\left( \frac{|x(s)|}{\psi(x(t-r))} \right) x(s), u(s) \right) ds$$

which directly implies

$$x(t-r+\delta) = x(t-r) + \int_{t-\delta}^{t} f\left( q\left( \frac{|x(s-r+\delta)|}{\psi(x(t-r))} \right) x(s-r+\delta), u(s-r+\delta) \right) ds \tag{A.16}$$

for all $t \ge \tilde{T}$. Exploiting (A.13), (A.14), (A.15), (A.10), (A.16) in conjunction with (4.9) we get for all $t \ge \tilde{T}$:

$$\begin{aligned} |\xi(t) - x(t-r+\delta)| &\le |z(t) - x(t-r)| + G_1 \int_{t-\delta}^{t} |\xi(s) - x(s-r+\delta)| ds + G_2 \int_{t-\delta}^{t} |z(s) - x(s-r)| ds \\ &+ \exp(-\mu t)\left| \xi(0) - z(0) - \int_{-\delta}^{0} f\left( q\left( \frac{|\xi(s)|}{\psi(z(0))} \right)\xi(s), u(s-r+\delta) \right) ds \right| \end{aligned}$$

The above inequality in conjunction with $G_1 \frac{\exp(\sigma\delta)-1}{\sigma} < 1$ implies that the following inequality holds for all $t \geq \tilde{T}$ :

$$
\begin{aligned}
&\sup_{\tilde{T} \leq s \leq t} \left( \exp(\sigma s) \left| \xi(s) - x(s-r+\delta) \right| \right) \leq \\
&\frac{G_1 \left( \exp(\sigma\delta) - 1 \right)}{\sigma - G_1 \left( \exp(\sigma\delta) - 1 \right)} \sup_{\tilde{T}-\delta \leq s \leq \tilde{T}} \left( \exp(\sigma s) \left| \xi(s) - x(s-r+\delta) \right| \right) \\
&+ \frac{\sigma + G_2 \left( \exp(\sigma\delta) - 1 \right)}{\sigma - G_1 \left( \exp(\sigma\delta) - 1 \right)} \sup_{\tilde{T}-\delta \leq s \leq t} \left( \exp(\sigma s) \left| z(s) - x(s-r) \right| \right) \\
&+ \frac{\sigma}{\sigma - G_1 \left( \exp(\sigma\delta) - 1 \right)} \left| \xi(0) - z(0) - \int_{-\delta}^{0} f\left( q\left( \frac{|\xi(s)|}{\psi(z(0))} \right) \xi(s), u(s-r+\delta) \right) ds \right|
\end{aligned}
\tag{A.17}
$$

Using (A.17), (4.2), (4.4) and (A.11) we are in a position to construct a non-decreasing function $M : \Re_+ \to \Re_+$ such that (4.11) holds. The proof is complete. ◁